\documentclass[ejs]{imsart}
  \setattribute{journal}{name}{}
  \setattribute{journal}{issn}{}

\usepackage{algorithm,algpseudocode}
\usepackage{amsthm,amsmath}
\RequirePackage[numbers]{natbib}
\RequirePackage[colorlinks,citecolor=blue,urlcolor=blue]{hyperref}
\usepackage{chngcntr}
\usepackage{graphicx,subcaption}
\usepackage{microtype}
\usepackage{yk}

\startlocaldefs

\def\exp#1{\text{\normalfont exp}\big( #1 \big)}

\def\b{\big}
\def\bb{\Big}
\def\bbb{\bigg}
\def\tX{\tilde{X}}
\def\bmu{\boldsymbol{\mu}}
\def\bzero{\mathbf{0}}
\def\ve{\varepsilon}
\endlocaldefs

\begin{document}

\begin{frontmatter}

\title{Statistical Convergence of the EM Algorithm on Gaussian Mixture Models}
\runtitle{Convergence of EM on Gaussian mixtures}


\begin{aug}
\author{\fnms{Ruofei} \snm{Zhao}}
\and
\author{\fnms{Yuanzhi} \snm{Li}}
\and
\author{\fnms{Yuekai} \snm{Sun}}

\address{Department of Statistics\\
University of Michigan \\
Ann Arbor, MI}

\runauthor{Zhao, Li, and Sun}

\affiliation{University of Michigan}

\end{aug}

\begin{abstract}
We study the convergence behavior of the Expectation Maximization (EM) algorithm on Gaussian mixture models with an arbitrary number of mixture components and mixing weights. We show that as long as the means of the components are separated by at least $\Omega(\sqrt{\min\{M,d\}})$, where $M$ is the number of components and $d$ is the dimension, the EM algorithm converges locally to the global optimum of the log-likelihood. Further, we show that the convergence rate is linear and characterize the size of the basin of attraction to the global optimum.
\end{abstract}





\end{frontmatter}

\section{Introduction}

The EM algorithm \citep{em} is an instrumental tool for evaluating the maximum likelihood estimator of latent variable models. It is a Majorization-Minimization (MM) algorithm that minimizes a surrogate objective function to avoid evaluating the intractable marginal (negative) log-likelihood of the latent variable model. However, as is shown by \citet{wu1983} and \citet{tseng}, the EM algorithm may not converge to a global minimizer of the log-likelihood. Instead, it may converge to a local minimizer or a stationary point. For a method that aims to evaluate the MLE, this is somewhat disappointing.

There is an recent line of research that shows the EM algorithm initialized in a neighborhood of the data generating parameters converges to the global minimizer.
Unfortunately, this line of work does not encompass the EM algorithm for fitting mixtures of more than two Gaussians. This paper fills this gap in the literature by providing conditions under which the EM algorithm for fitting Gaussian mixture models with an arbitrary number of well-separated  components and arbitrary mixing weights converges to the global minimizer. We show that the EM algorithm converges linearly as long as it is initialized in a neighborhood of the true centers. Our results are of the same flavor as those by \citet{bowei} and \citet{balakrishnan2017}, and our proofs follow the same general route.

This paper is organized as follows. The rest of this section briefly reviews related work on the EM algorithm. Section 2 describes the EM algorithm for fitting mixtures of Gaussians and introduces a population version of the algorithm that appears in our study. Section 3 states our main results on the convergence of EM. In section 4, we present simulation results that validate some of our theoretical results. In section 5, we prove the main results. Finally, in section 6, we discuss our results and compare them to similar results in the literature.

\subsection{Related work}


Most closely related to our work is the line of work on the convergence of the EM algorithm for fitting Gaussian mixture models (GMM). On a mixture of two equally-weighted Gaussians, \citet{balakrishnan2017} first derived statistical convergence results by specializing the general framework they proposed to this model. Their framework  also applies to a variant of the EM algorithm known as gradient EM, and they used their framework to obtain similar results for gradient EM. \citet{brinda16} later obtained results of a similar flavor, but showed that there is a larger neighborhood of the true centers within which the EM algorithm converges. Finally, \citet{xuji} and \citet{daskalakis} completely characterized the global convergence behavior of the EM algoritm for fitting two equally weighted Gaussians. When there are more than two  components, a result by \citet{badmaxima} showed that bad local minima exists even in the idealized case of equally weighted mixtures of well-separated spherical Gaussians. \citet{bowei} proved local convergence results for the gradient EM algorithm for fitting mixtures of an arbitrary number of well-separated spherical Gaussians. Despite the recent progress, we are not aware of any results that characterize the local convergence behavior of the EM algorithm on mixtures of two or more Gaussians. For variants of the EM algorithm for fitting high-dimensional mixture models, we refer readers to \citet{dasgupta07,chime,highdem,regem}.

There is also a large body of work on other methods for learning mixtures of Gaussians and, more generally, finite mixtures. One major line of work \citep{dasgupta99,vempala04,achlioptas,arora2005,vempala08,Brubaker08} is based on  dimension reduction techniques (such as spectral embeddings). Like the EM algorithm, these methods require the centers of the mixture components to be well-separated. Another more recent line of work \citep{Belkin10,Kalai10,Moitra10,Hsu13,Hardt15} employs the method-of-moments, and allows the centers of the mixture components to be arbitrarily close (as long as the sample size is large enough). Some other important algorithms and theoretical work are \citet{Brubaker08,Chaudhuri08,Chaudhuri09,yulu}. For statistical properties such as the rate of convergence of the MLE or the rate of convergence of the estimated mixing distribution, we refer readers to \citet{ghosal2001, nguyen2013,heinrich2018} and the references therein.

\section{The EM algorithm on Gaussian mixture models}

We consider Gaussian mixture models with known mixture weights and known common covariance structure. Formally, suppose there are $M$ isotropic Gaussian distributions, $\cN(\bmu_1^*,I_d),\dots,\cN(\bmu_M^*,I_d)$, and mixture weights, $\pi_1,\dots,\pi_M \ge 0$. The Gaussian mixture model we consider is the set of densities 
\begin{equation}
\{{\textstyle p(x;\bmu^*) = \sum_{i=1}^M\pi_i\phi(x - \bmu_i^*):\bmu^* = \big[(\bmu_1^*)^T, \dots, (\bmu_M^*)^T\big]^T\in\reals^{Md}}\},
\label{eq:gaussianMixtureModel}
\end{equation}
where $\phi(z) = \frac{1}{(2\pi)^{1/2}}e^{-\frac12\|z\|_2^2}$ is the standard Gaussian density in $\reals^d$. The assumption that the components are isotropic leads to no loss of generality as long as (i) the mixture components share a common covariance structure and (ii) this structure is known. Without loss of generality, we also assume mixture is centered; \ie\ it has mean zero. The task of fitting \eqref{eq:gaussianMixtureModel} boils down to estimating $\bmu^*$ from observations
\[
X_1,\dots,X_n\overset{\iid}{\sim}p(\cdot;\bmu^*).
\]

The {\it EM algorithm} for fitting a Gaussian mixture model alternates between evaluating the posterior probabilities (E-step) of the labels and updating the estimates of the parameters (M-step). We combine the two steps to arrive at the update rule
\begin{equation}
\bmu^+_i \gets \frac{\sum_{j=1}^{n} w_i(X_j;\bmu)X_j}{\sum_{j=1}^{n} w_i(X_j;\bmu)},\quad i = 1,\ldots,M,
\label{eq:updaterule1}
\end{equation}
where the weights $w_i(X_j;\bmu)$ are defined as
\begin{equation}
w_i(x;\bmu) := \frac{\pi_i \phi(x;\bmu_i)}{\sum_{k=1}^M \pi_k \phi(x;\bmu_k) }.
\label{eq:definew}
\end{equation}
We see that $w_i(X;\bmu)$ is the probability that $X$ comes from component $i$ and $\bmu_i^+$ is a weighted average of samples, where the weights are the $w_i(X_j;\bmu)$'s. For this reason, the EM algorithm is known as a soft version of the $K$-means algorithm.

In our analysis, we work with a population version of the EM algorithm. Its update rule is:
\begin{equation}
\bmu^+_i \gets \frac{\bbE[ w_i(X;\bmu)X]}{\bbE w_i(X;\bmu)} = \frac{\int_{\reals^d}w_i(x;\bmu)xp(x;\bmu^*)dx}{\int_{\reals^d}w_i(x;\bmu)p(x;\bmu^*)dx},\quad i = 1,\ldots,M.
\label{eq:updaterule2}
\end{equation}
We emphasize that the expectations in \eqref{eq:updaterule2} are with respect to the true data generating process. We will first derive convergence results for this population version of the EM algorithm and then extend this result to (the sample-version of) the EM algorithm using concentration results.

{\bf Notation:} We summarize the notations in this paper: $d$ is the dimension of the observations $X_1,\dots,X_n$ and $M$ is the number of mixture components. Define $R_{\max}$, $R_{\min}$ as the largest and smallest distances between the centers of any pair of mixture components: 
\[
\begin{aligned}
R_{\max}=\max_{i \neq j} \|\bmu^*_i - \bmu^*_j \|_2, \\
R_{\min}=\min_{i \neq j} \|\bmu^*_i - \bmu^*_j \|_2.
\end{aligned}
\]
Define $\kappa$ as the smallest mixture weight: $\kappa = \min\{\pi_1,\dots,\pi_M\}$. Given two positive sequences $\{a_n\},\{b_n\}$, $a_n = \cO(b_n)$ means there exists an absolute constant $C$ such that $a_n \leq C b_n$ for all $n$; $a_n = \Omega(b_n)$ there exists an absolute constant $C$ such that $a_n \geq C b_n$ for all $n$. We write $a_n = \Theta(b_n)$ if $a_n = \cO(b_n)$ and $a_n = \Omega(b_n)$; we write $a_n = o(b_n)$ if $\frac{a_n}{b_n} \rightarrow 0$ as $n \rightarrow \infty$.

\section{Statement of the main results}
\label{sec:mainResutls}

In this section, we state our main results for the convergence of the EM algorithm on Gaussian mixture models. 
First, we show that the population version of the EM algorithm converges linearly to $\bmu^*$ as long as (i) certain signal strength conditions are met and (ii) the algorithm is initialized in a neighborhood of $\bmu^*$. We also characterize the size of this neighborhood in terms of the properties of the data generating process.

\begin{theorem}
 Suppose $R_{\min} \geq 30\min\{2M,d\}^{\frac12}$ and radius $a$ satisfies
  \begin{equation}
   a \leq \frac12R_{\min} - \min\{d,2M\}^{\frac12}\max\{{\textstyle 4\sqrt{2}[\log(\frac{R_{\min}}{4})]_+^{\frac12},8\sqrt{3},8\log (\frac{4}{\kappa})}\}.
  \label{eq:contractionradius}
  \end{equation}
  Then for any iterate $\bmu$ satisfying $\max_{i \in [M]}\|\bmu_i-\bmu^*_i\|_2 \leq a$, the next iterate $\bmu^+$ given by \eqref{eq:updaterule2} satisfies
  \[
  \max_{i \in [M]} \|\bmu^+_i-\bmu^*_i\|_2 \leq \zeta \max_{i \in [M]} \|\bmu_i-\bmu^*_i\|_2,
  \]
  where 
  \begin{equation}
  \zeta = \frac{3M}{\kappa^2}( 2 R_{\max} + \min\{2M,d\})^2\text{\normalfont exp}({\textstyle - \frac18(\frac12R_{\min} -a)\min\{d,2M\}^{\frac12}}).
  \label{eq:zeta}
  \end{equation} 
  \label{thm:populationEMConvergence}
\end{theorem}

We remark that Theorem \ref{thm:populationEMConvergence} does not imply cluster-wise convergence; \ie\ Theorem \ref{thm:populationEMConvergence} does not imply
\[
\|\bmu_i^+ - \bmu_i^*\|_2 \le \zeta \|\bmu_i-\bmu^*_i\|_2\text{ for all }i\in[M].
\]
This is hardly surprising. If we initialize the EM algorithm at 
\[
\big[(\bmu_1^*)^T,\bmu_2^T,\ldots,\bmu_M^T\big]^T,
\]
where $\bmu_2,\dots,\bmu_M$ are arbitrary points in $\reals^d$, there is no guarantee that $\bmu_1^+$ remains at $\bmu_1^*$. 

Recall the mixture components are isotropic, so $R_{\min}$ is the signal strength. Theorem \ref{thm:populationEMConvergence} requires the signal strength to grow as $\Omega(\min\{M,d\}^{\frac12})$. Theorem \ref{thm:populationEMConvergence} also shows that the contraction coefficient $\zeta$ is decreasing in $R_{\min}$ and $\kappa$. It also shows that the contraction radius $a$ approaches $\frac12R_{\min}$ as $R_{\min}$ goes to infinity. 
This contraction radius is essentially optimal because there are examples of the EM algorithm converging to non-global local minima if $\|\bmu_i-\bmu^*_i\|=\frac12R_{\min}$. 
For example, consider the task of fitting a mixture of two Gaussians. If we initialize the EM algorithm at 
\[\textstyle
\bigl[\frac12(\bmu^*_1 + \bmu^*_2)^T,\frac12(\bmu^*_1 + \bmu^*_2)^T\big]^T,
\]
the algorithm will never separate the centers and converges to a stationary point with identical estimates of the two centers. It is not clear whether the $\min\{d,2M\}^{\frac12}$ term in \eqref{eq:contractionradius}{} is optimal. We suspect not, because in the experiments we run, we have never seen EM fail to converge to the truth when the initializer satisfies $\max \|\bmu_i-\bmu^*_i \|_2 < R_{\min}/2$.
The improved dependence of the contraction radius on $\min\{d,2M\}$ instead of $d$ is intuitive. Informally, if $d > 2M$, the components of the noise orthogonal to the subspace spanned by the centers and estimated centers cancel out when expectation is taken, which reduces the effective dimension of the problem from $d$ to $\min\{d,2M\}$. The details are in the the proof of Theorem \ref{thm:populationEMConvergence}. We remark that this improvement is not present in studies of the convergence of the EM algorithm that does not go through a population level analysis (\cf\ the result of \citet{dasgupta07}).

In the statement of theorem \ref{thm:populationEMConvergence}, there is no provision that $\zeta < 1$. Rearranging $\zeta < 1$ for conditions on $a$ leads to Theorem \ref{cor:populationEMConvergence}, which is easier to parse but obscures the dependence of the contraction coefficient on the properties of the data generating process.

\begin{theorem}
As long as $R_{\min} \ge C_0\min\{d,M\}^{\frac12}$ and
\[
\max_{i \in [M]}\|\bmu_i-\bmu^*_i\|_2 \leq \frac12R_{\min} - C_1\min\{d,M\}^{\frac12}\log( \max\{{\textstyle\frac{M}{\kappa^2},  R_{\max}, \min\{d,M\} }\})^{\frac12},
\]
where $C_0,C_1 > 0$ are universal constants, we have
  \begin{equation}
  \max_{i \in [M]} \|\bmu^+_i-\bmu^*_i\|_2 \leq \frac12 \max_{i \in [M]} \|\bmu_i-\bmu^*_i\|_2.
  \label{eq:thm1,2}
  \end{equation}
  \label{cor:populationEMConvergence}
\end{theorem}

Second, we carry out a perturbation analysis to extend Theorem \ref{cor:populationEMConvergence} to (the sample version of) the EM algorithm. At a high level, our result states that the iterates of the EM algorithm converge linearly to $\bmu^*$ up to the statistical precision of the estimation task. Our proof hinges on recent concentration results by \citet{songmei}. Compared to the proof of an analogous result for the gradient EM algorithm by \citet{bowei}, our proof is considerably simpler.

\begin{theorem}
Suppose $R_{\min} \geq C_0{\min\{d,M\}}^{\frac12}$ and the initial iterate $\bmu^0$ satisfies
  \begin{equation}
  \max_{i \in [M]} \|\bmu^0_i-\bmu^*_i\|_2 \leq \frac12R_{\min} - C_1\min\{d,M\}^{\frac12}\log( \max\{{\textstyle\frac{M}{\kappa^2},  R_{\max}, \min\{d,M\} }\})^{\frac12},
  \label{eq:sampleradius}
  \end{equation}
  where $C_0, C_1 > 0$ are universal constants. As long as the sample size $n$ is large enough so that
  \begin{equation}
  \frac{\log n}{n} \leq \min \bbb\{\frac{\kappa^2}{144\tilde{C_2}Md},\frac{\kappa^2\max_{i \in [M]} \|\bmu^0_i-\bmu^*_i\|_2^2}{9\tilde{C_3} R_{\max}^2 Md}\bbb\},
  \label{eq:samplesizecondition}
  \end{equation}
  where $C_2,C_3 > 0$ are universal constants and $\tilde{C_2}= C_2\log (M(2R_{\max}+\sqrt{d}))$, $\tilde{C_3}=C_3\log (M(3R_{\max}^2+d))$, the subsequent EM iterates $\{\bmu^t \}_{t=1}^\infty$ given by \eqref{eq:updaterule1} satisfy
  \begin{equation}
  \max_{i \in [M]} \|\bmu^t_i-\bmu^*_i\|_2 \leq \frac{1}{2^t} \max_{i \in [M]} \|\bmu^0_i-\bmu^*_i\|_2+ {\textstyle\frac{3R_{\max}}{\kappa}(\frac{\tilde{C_3}Md\log n}{n}})^{\frac12}
  \label{eq:samplecontraction}
  \end{equation}
  with probability at least $1-\frac{2M}{n}$.
  \label{thm:EMConvergence}
\end{theorem} 

We see that the first term on the right side of \eqref{eq:samplecontraction} converges to zero linearly while the second term does not depend on $t$. Initially, the first term on the right side dominates the second term and the right side decreases linearly. However, after sufficiently many iterations, the second term on the right side dominates the first term, and the right side settles down to a limit that is $\cO((\frac{Md}{\kappa^2 n})^{\frac12})$. We recognize this limit as the statistical precision (modulo constant and logarithmic factors).


\textbf{Convergence to the maximum likelihood estimator.} We remark that Theorem \ref{thm:EMConvergence} implies EM converges to a stationary point within a ball of radius $\cO((\frac{Md}{\kappa^2n})^{\frac12})$ of the true centers. It is known that the maximum likelihood estimator (MLE) falls inside this ball with high probability. As long as there are no other stationary points in this ball, Theorem \ref{thm:EMConvergence} implies EM converges to the MLE. This is a consequence of the gradient stability result by \citet{bowei}, which implies $\bmu^*$ is the only stationary point of the expected log-likelihood function in $\otimes_i (\bmu^*_i, a)$, and concentration results by \citet{songmei}, which imply the log-likelihood function only has one stationary point in $\otimes_i (\bmu^*_i, a)$.

\section{Simulation results}

\begin{figure*}[!htbp]
  \captionsetup{width=\linewidth}
  \centering
  \begin{subfigure}[b]{0.46\textwidth}
    \includegraphics[width=\textwidth]{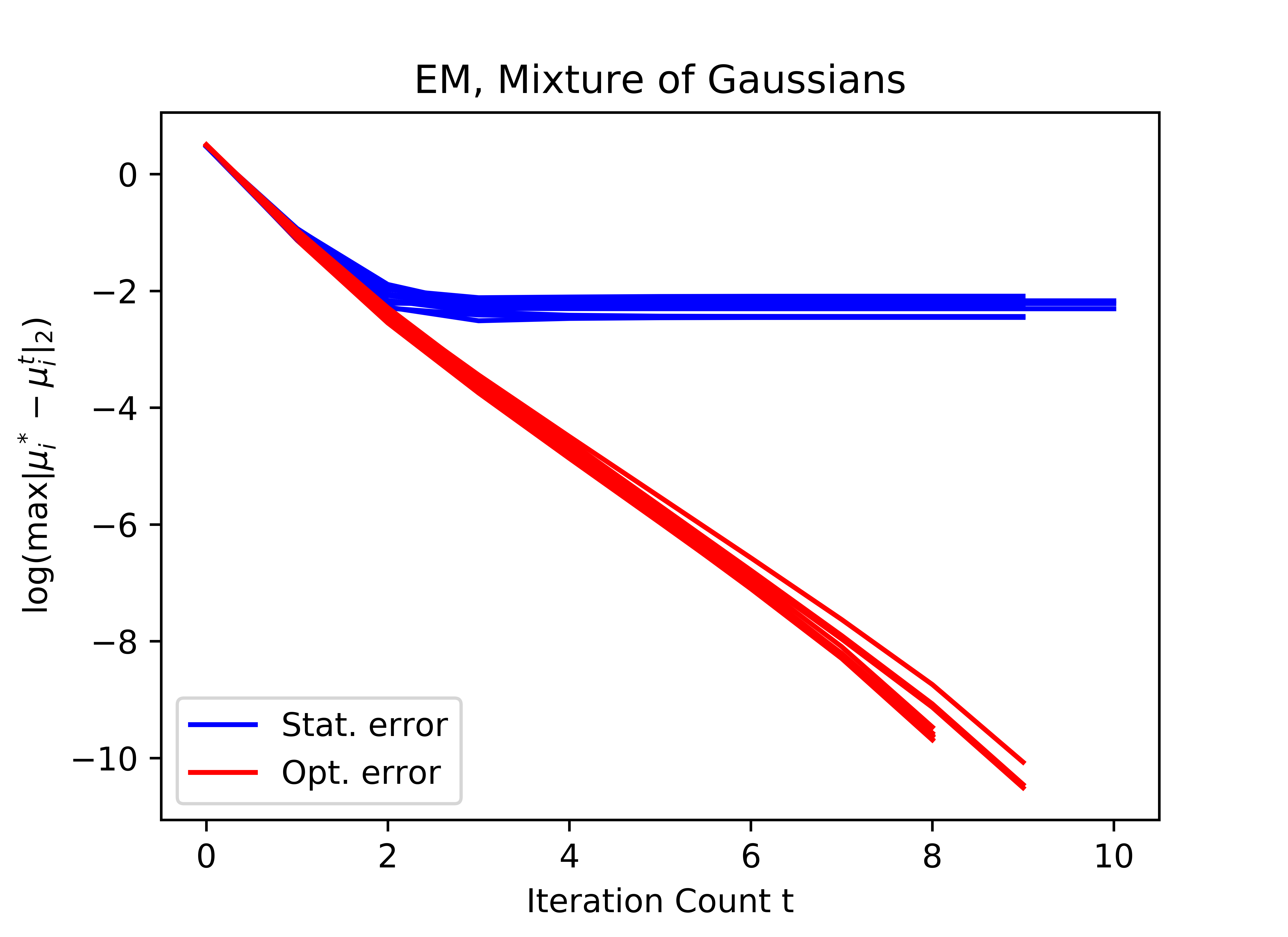}
    \caption{Figure 1.a}
  \end{subfigure}
  ~ \qquad 
  \begin{subfigure}[b]{0.46\textwidth}
    \includegraphics[width=\textwidth]{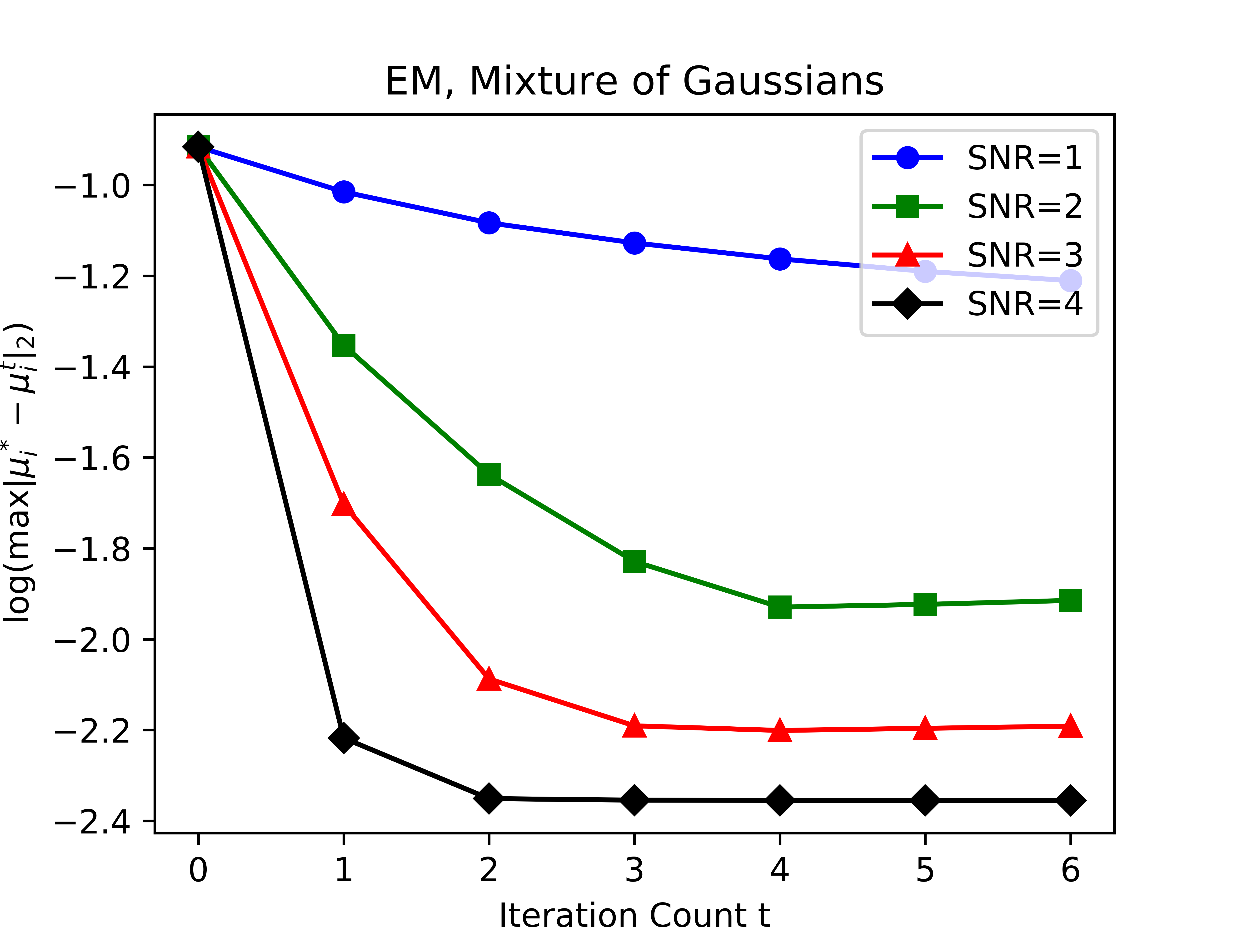}
    \caption{Figure 1.b}
  \end{subfigure}\\
\begin{subfigure}[b]{0.46\textwidth}
  \includegraphics[width=\textwidth]{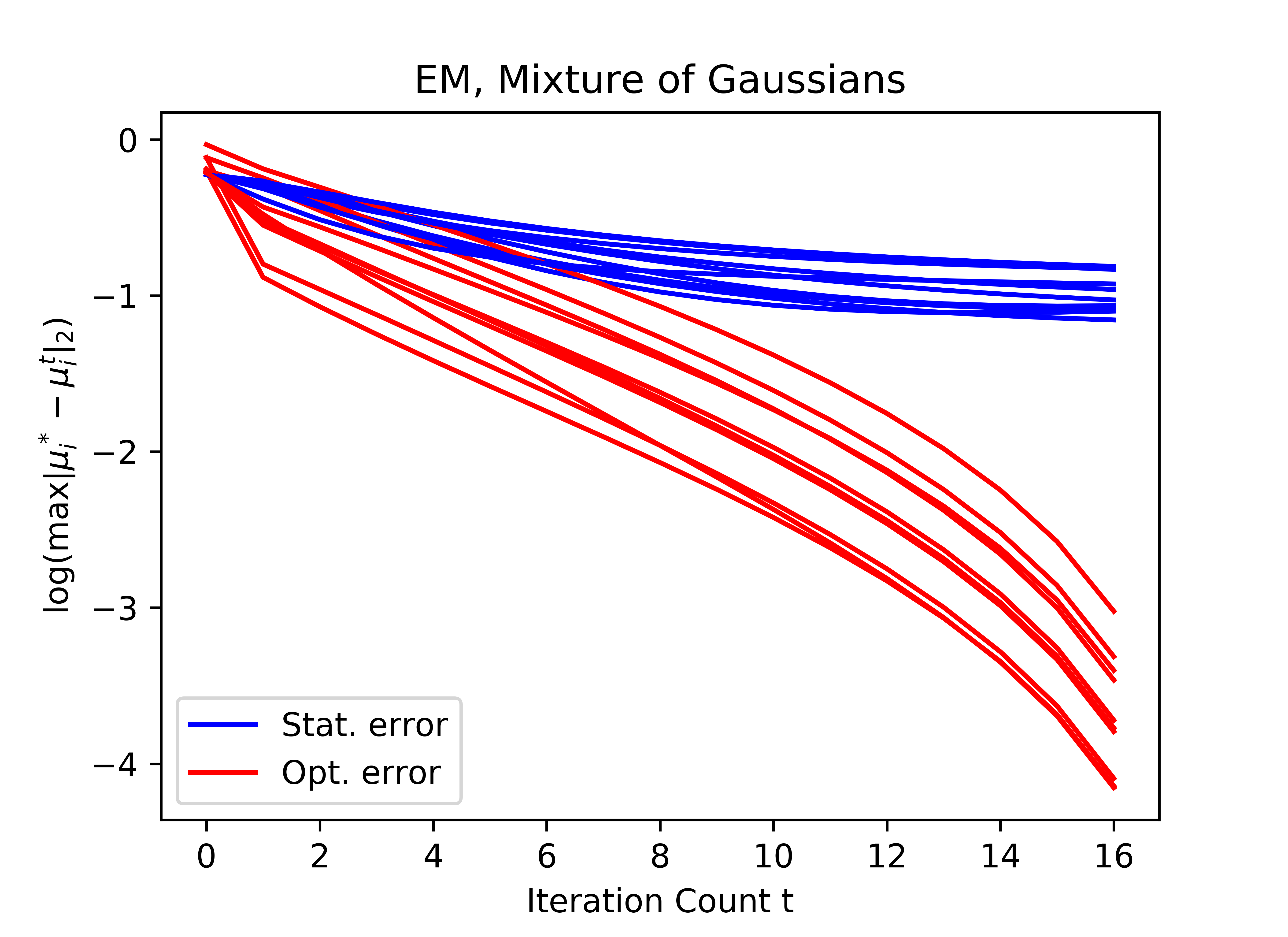}
  \caption{Figure 1.c}
\end{subfigure}
~ \qquad 
\begin{subfigure}[b]{0.46\textwidth}
  \includegraphics[width=\textwidth]{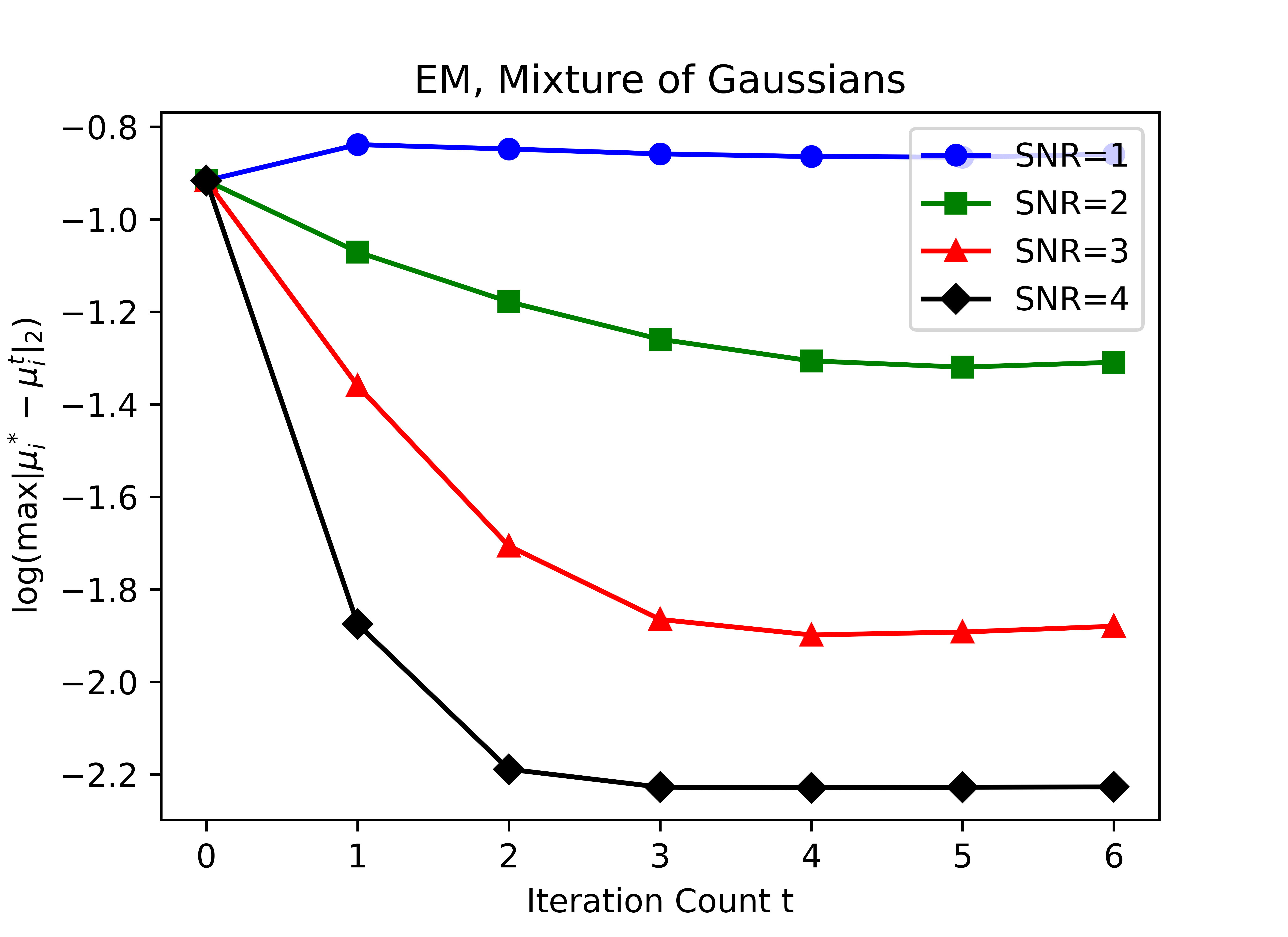}
  \caption{Figure 1.d}
\end{subfigure}

  \caption{Convergence of EM on the mixture of $M=5$ Gaussians in $\bbR^{10}$. Left column shows the statistical error and optimization error for ten trials. Right column shows average statistical error over ten trials under different SNR (SNR is $R_{min}$). While top row has balanced cluster weights, bottom row does not.}\label{fig:experiments}
\end{figure*}

In this section, we present some simulation results to validate our theoretical results. In the first set of experiments, we check the prediction of Theorem \ref{thm:populationEMConvergence} that the statistical error decreases linearly initially and eventually reaches a plateau. The data generating process is a mixture of $M = 5$ isotropic Gaussians in $\bbR^{10}$ with one mixture component centered at the origin and the remaining four centers at the vertices of $R_{\min}\Delta^9$, where $\Delta^9$ is the probability simplex in $\reals^{10}$, so $\frac{R_{\max}}{R_{\min}} = \sqrt{2}$. The components are equally weighted ($\pi_i = \frac15$). 

The top left panel of Figure \ref{fig:experiments} shows the decrease of the optimization error $\max_i \|\hat{\bmu}_i - \bmu_i^t\|_2$ and statistical error $\max_i \|\bmu^*_i - \bmu_i^t\|_2$ over 10 runs of the EM algorithm. We set $R_{\min} = 2$ and generate $n = 8000$ samples from the Gaussian mixture model. We initialize the EM algorithm at 
\[
\big[(\bmu_1^* + \delta_1)^T,\dots,(\bmu_M^* + \delta_M)^T\big]^T,
\]
where $\delta_i$ is uniformly distributed on the sphere of radius $0.4\cdot R_{\min}$. We see that the statistical error decreases linearly initially and eventually reaches a plateau after four iterations. This agrees with the implications of Theorem \ref{thm:EMConvergence}. The top right panel of Figure \ref{fig:experiments} shows the average log statistical errors over 10 trails as $R_{\min}$ varies. We see that larger $R_{\min}$ values lead to faster convergence, which agrees with the implications of Theorem \ref{thm:populationEMConvergence}.

The plots in the bottom panels are analogous to the plots above. We keep all simulation parameters, except we change the mixture weights from uniform to $\pi_i = \frac{i}{15}$, $i\in[M]$. We see that non-uniform mixture weights hurt the performance of the EM algorithm. Comparing the two plots on the left, we see that non-uniform weights causes the algorithm to converge slower and reduces the statistical precision of the output. The slower convergence agrees with the implications of Theorem \ref{thm:populationEMConvergence}, which shows that the contraction coefficient is inversely proportional to the minimum mixture weight $\kappa$. The reduced statistical precision is due to the fact that the centers of mixture components with smaller weights are estimated less accurately. This is because the mixture components with smaller weights have smaller effective sample sizes. We also see greater variation across the ten runs of the algorithm.

\section{Proofs of the main results}

We prove Theorem \ref{thm:populationEMConvergence}, \ref{cor:populationEMConvergence}, and \ref{thm:EMConvergence} in this section, deferring the proofs of the technical lemmas to the appendices.

\subsection{Proof of Theorem \ref{thm:populationEMConvergence}}

We make a few observations before proceeding to the proof. Without loss of generality, we focus on the update rule for the first center $\bmu_1$:
\begin{equation*}
\bmu^+_1-\bmu^*_1 = \frac{\bbE\big[ w_1(X;\bmu)(X-\bmu^*_1)\big]}{\bbE\bigl[w_1(X;\bmu)\bigr]}
\end{equation*}
The vector of true centers $\bmu^*$ is a fixed point of \eqref{eq:updaterule2}, which implies
\[
\bbE\big[ w_1(X;\bmu^*)(X-\bmu^*_1)\big] = \bzero.
\]
Thus
\begin{equation}
\bmu^+_1-\bmu^*_1 = \frac{\bbE\big[ \big(w_1(X;\bmu) - w_1(X;\bmu^*) \big)(X-\bmu^*_1)\big]}{\bbE\big[w_1(X;\bmu)\big]}.
\label{eq:dif}
\end{equation}
In the first step of the proof, we establish an upper bound on the norm of the numerator. In the second step, we establish a lower bound on the denominator in \eqref{eq:dif}. Finally, in the third step, we combine the upper and lower bounds to show the EM update rule is a contraction.

\textbf{Step 1 (upper bounding the numerator).} Define $\bmu^t := \bmu^* + t(\bmu-\bmu^*)$ and define $g_X(t) := w_1(X;\bmu^t)$. We have
\[
\begin{aligned}
w_1(X;\bmu)-w_1(X;\bmu^*) &= \int_0^1 g_X'(t) dt = \int_0^1  \nabla_{\bmu} w_1\b(x;\bmu^t\b)^T(\bmu^t-\bmu^*)dt,
\end{aligned}
\]
where
{
  \begin{equation}
  \nabla_{\bmu} w_1(X;\bmu) =\\
  \begin{bmatrix}
  -w_1(X;\bmu)(1-w_1(X;\bmu))(\bmu_1-X)\\
  w_1(X;\bmu)w_2(X;\bmu)(\bmu_2-X)\\
  \vdots\\
  w_1(X;\bmu)w_M(X;\bmu)(\bmu_M-X)\\
  \end{bmatrix}.
  \end{equation}
}
We thus have 
\begin{align}
& \big\| \bbE\big[(w_1(X;\bmu) - w_1(X;\bmu^*))(X-\bmu^*_1)\big] \big\|_2 \nonumber \\
&\quad\leq \bb\| \int_{0}^1 \bbE\big[w_1(X;\bmu^t)(1-w_1(X;\bmu^t))(X-\bmu^t_1)^T(\bmu_1^t-\bmu_1^*)(X-\bmu_1^*)\big]dt \nonumber\\
&\qquad - \sum_{i\neq 1}{\textstyle\int_{0}^1 \bbE\big[w_1(X;\bmu^t)w_i(X;\bmu^t)(X-\bmu^t_i)^T(\bmu_i^t-\bmu_i^*)(X-\bmu_1^*)\big]dt} \bb\|_2 \nonumber\\
&\quad\leq \int_{0}^1 \b\| \bbE\big[w_1(X;\bmu^t)(1-w_1(X;\bmu^t))(X-\bmu_1^*) (X-\bmu^t_1)^T\big] \b\|_{op} \| \bmu_1^t-\bmu_1^* \|_2dt\nonumber\\
&\qquad + \sum_{i\neq 1}{\textstyle\int_{0}^1 \b\| \bbE\big[w_1(X;\bmu^t)w_2(X;\bmu^t)(X-\bmu_1^*)(X-\bmu^t_i)^T\big] \b\|_{op}  \| \bmu_i^t-\bmu_i^* \|_2dt} \nonumber\\
&\quad\leq V_1 \| \bmu_1-\bmu_1^* \|_2 + \sum_{i\neq 1} V_i \| \bmu_i-\bmu_i^* \|_2 \nonumber\\
&\quad\leq M ({\textstyle\max_i V_i}) \cdot \max_i \| \bmu_i-\bmu_i^* \|_2, \label{eq:integrationexpression}
\end{align}
where
\begin{align}
  V_1 =\sup_{t \in [0,1]} \b\| \bbE\big[w_1(X;\bmu^t)(1-w_1(X;\bmu^t))  (X-\bmu_1^*)  (X-\bmu^t_1)^T\big] \b\|_{op}, \label{eq:V1}\\
  V_i = \sup_{t \in [0,1]} \b\| \bbE\big[w_1(X;\bmu^t)w_i(X;\bmu^t)(X-\bmu_1^*)  (X-\bmu^t_i)^T\big]\b\|_{op}. \label{eq:V2}
\end{align}
The rest of the first step consists of establishing bounds on the $V_i$'s. We state the result here and defer the details to Appendix A.

\begin{lemma}
  As long as $R_{\min} \geq 30\min\{d,2M\}^{\frac12}$ and 
  \begin{equation}
  a \leq \frac12R_{\min} - \min\{d,2M\}^{\frac12} \cdot \max\{4\sqrt{2}[\log(R_{\min}/4)]_+^{\frac12},8\sqrt{3}\},
  \label{eq:conditionfora_1}
  \end{equation}
  then for any $\bmu$ such that $\bmu_i \in \cB(\bmu_i^*, a), i \in [M]$, we have 
  \begin{equation}
  \begin{aligned}
  \max_{i \in [M]} V_i \leq \frac{2}{\kappa} ( 2 R_{\max} + \min\{2M,d\})^2  \text{\normalfont exp}({\textstyle-\frac18( \frac12R_{\min} -a)  \min\{d,2M\}^{\frac12}}).
  \end{aligned}
  \label{eq:boundmaxv}
  \end{equation}
  \label{lemma:upperbound}
\end{lemma}

\textbf{Step 2 (lower bounding the denominator).} We state the lower bound and describe the underlying intuition, deferring the proof to the appendix. Let $Z$ be the label of $X$. We observe that
\[
\bbE\big[w_1(X;\bmu^*)\big] = \bbE\big[\bbP_{\bmu^*}(Z=1|X)\big]=\pi_1 > \kappa.
\]
As long as $\bmu\approx \bmu^*$, $\Ex\big[w_1(X;\bmu)\big]\approx\Ex\big[w_1(X;\bmu^*)\big]$, so $\bbE\big[w_1(X;\bmu)\big]$ cannot be much smaller than $\kappa$. 

\begin{lemma}
As long as $R_{\min} \geq 30\min\{M,d\}^{\frac12}$ and 
  \begin{equation}
  a \leq \frac12R_{\min} - \min\{M,d\}^{\frac12} \cdot \max\{{\textstyle 4\sqrt{2}[\log(\frac{R_{\min}}{4})]_+^{\frac12},8\sqrt{3}, 8\log(\frac{4}{\kappa})}\},
  \label{eq:conditionfora_2}
  \end{equation}
  then for any $\bmu$ such that $\bmu_i \in \cB(\bmu_i^*, a), i \in [M]$, we have 
  \begin{equation}
  \bbE\big[w_i(X;\bmu)\big] \geq \frac{3}{4} \kappa,\quad i\in[M].
  \end{equation}
  \label{lemma:lowerbound}
\end{lemma} 

\textbf{Step 3:} We combine the bounds for $\max_i V_i$ and $\bbE\big[w_1(X;\bmu)\big]$ to show that the EM update rule is a contraction. Without loss of generality, we focus on the update rule for the first cluster. By \eqref{eq:dif}, \eqref{eq:integrationexpression}, and Lemma \ref{lemma:lowerbound}, we have
{
  \begin{equation}
  \begin{aligned}
  \| \bmu^+_1-\bmu^*_1 \|_2&= \frac{\| \bbE[ \big(w_1(X;\bmu) - w_1(X;\bmu^*) \big)(X-\bmu^*_1)] \|_2}{\bbE\big[w_1(X;\bmu)\big]}\\
  & \leq \frac{4M}{3\kappa}({\textstyle\max_i V_i})({\textstyle\max_i \| \bmu_i-\bmu_i^* \|_2}).
  \end{aligned}
  \label{eq:contractionpop}
  \end{equation}}
Plugging in \eqref{eq:boundmaxv}, we have
\begin{equation}
\begin{aligned}
\| \bmu^+_1-\bmu^*_1 \|_2 \leq &\frac{8M}{3\kappa^2}( 2 R_{\max} + \min\{2M,d\})^2 \\
&\cdot \text{\normalfont exp}({\textstyle -\frac18(\frac12R_{\min} -a )  \min\{d,2M\}^{\frac12}}) (\max_i \| \bmu_i-\bmu_i^* \|_2).
\end{aligned}
\end{equation}
We recognize the factor in front of $\max_i \| \bmu_i-\bmu_i^* \|_2$ is smaller than the contraction factor $\zeta$ in Theorem \ref{thm:populationEMConvergence}. We can also check that \eqref{eq:contractionradius} implies \eqref{eq:conditionfora_1} and \eqref{eq:conditionfora_2} so the conditions of Lemmas \ref{lemma:upperbound} and \ref{lemma:lowerbound} are satisfied. Putting the two parts together, we see that Theorem \ref{thm:populationEMConvergence} is correct.

\subsection{Proof of Theorem \ref{cor:populationEMConvergence}}

To prove Theorem \ref{cor:populationEMConvergence}, we start from \eqref{eq:contractionpop} and solve for the contraction radius $a$ that implies $\zeta < \frac12$. It is enough to find $\alpha$ such that 
\begin{equation*}
\frac{M \max_i V_i}{\frac23\kappa} \leq \frac{1}{2}.
\end{equation*}
By Lemma \ref{lemma:upperbound}, the above condition is implied by
\[
\frac{2}{\kappa}( 2 R_{\max} + \min\{2M,d\})^2 \text{\normalfont exp}({\textstyle -\frac18( \frac12R_{\min} -a)  \min\{d,2M\}^{\frac12}}) \leq \frac{\kappa}{3M}.
\]
Rearranging, we have
\begin{equation}
a \leq \frac12R_{\min} - \frac{8}{\min\{d,2M\}^{\frac12}}\log({\textstyle\frac{6M( 2 R_{\max} + {\min\{2M,d\}} )^2}{\kappa^2}}).
\label{eq:conditionfora_3}
\end{equation}
Finally, we check that there is a universal constant $C_1$ such that
\[
a  \leq \frac12R_{\min} - C_1\min\{d,2M\}^{\frac12}\log( \max\{{\textstyle\frac{M}{\kappa^2},  R_{\max}, \min\{2M,d\}}\})^{\frac12},
\]
implies \eqref{eq:conditionfora_1}, \eqref{eq:conditionfora_2}, and \eqref{eq:conditionfora_3}.

\subsection{Proof of Theorem \ref{thm:EMConvergence}}

The intuition underlying the proof of Theorem \ref{thm:EMConvergence} is the population and sample update rules (\eqref{eq:updaterule1} and \eqref{eq:updaterule2} respectively) are similar. In the proof, we appeal to the following technical lemmas on the uniform convergence of $\frac{1}{n} \sum_{j=1}^nw_i(X_j;\bmu)(X_j-\bmu^*_1)$ and $\frac{1}{n} \sum_{j=1}^nw_i(X_j;\bmu)$ to their population counterparts. 

\begin{lemma}
Define the event $\cE_{1,i}$ as 
\[
\begin{aligned}
&\textstyle\sup_{\bmu \in \cU}\| \frac{1}{n}\sum_{j=1}^nw_i(X_j;\bmu)(X_j-\bmu^*_i) - \bbE\big[w_i(X;\bmu)(X-\bmu^*_i)\big]\|_2 \\
&\textstyle\qquad\geq 1.5R_{\max}(\frac{\tilde{C_3}Md\log n}{n})^{\frac12},
\end{aligned}
\]
where $\cU=\otimes_{i=1}^M \cB(\bmu^*, 1.5R_{\max})$, $C_3$ is a universal constant, and 
\[
\tilde{C_3}=C_3\log (M(3R_{\max}^2+d)).
\]
We have $\bbP(\cE_{1,i}) \leq \frac{1}{n}$.
\label{lemma:contractionvector}
\end{lemma}

\begin{lemma}
Define the event
\[\textstyle
\cE_{2,i} = \sup_{\bmu \in \cU}| \frac{1}{n}\sum_{j=1}^n w_i(X_j,\bmu) - \bbE\big[w_i(X;\bmu)\big]| \geq  (\frac{\tilde{C_2}Md\log n}{n})^{\frac12},
\]
where $\bmu\in\cU=\otimes_{i=1}^M \cB(\bmu^*, R_{\max})$, $C_2$ is a universal constant, and 
\[
\tilde{C_2}=C_2\log (M(2R_{\max}+\sqrt{d})).
\]
We have $\bbP(\cE_{2,i}) \leq   \frac{1}{n}$.
  \label{lemma:contractionw}
\end{lemma}

\noindent \textbf{Proof of Theorem \ref{thm:EMConvergence}:}
On the event $\{\cap_{i\in[M]}\cE_{1,i}\}\cap\{\cap_{i\in[M]}\cE_{2,i}\}$, we have
\begin{align*}
&\textstyle\sup_{\bmu \in \cU} | \frac{1}{n} \sum_{j=1}^n w_i(X_j,\bmu) - \bbE\big[w_i(X;\bmu)\big]| \leq (\frac{\tilde{C_2}Md\log n}{n})^{\frac12},\\
&\textstyle\sup_{\bmu \in \cU} \| \frac{1}{n}\sum_{j=1}^n w_i(X_j;\bmu)(X_j-\bmu^*_i) - \bbE\big[w_i(X;\bmu)(X-\bmu^*_i)\big]\|_2 \\
&\quad\leq 1.5R_{\max}(\frac{\tilde{C_3}Md\log n}{n})^{\frac12}
\end{align*}
for all $i\in[M]$. By Lemmas \ref{lemma:contractionvector} and \ref{lemma:contractionw}, this event occurs with probability at least $1-\frac{2M}{n}$.  The minimum sample size condition \eqref{eq:samplesizecondition} implies
\begin{align} 
R_{\max}({\textstyle\frac{\tilde{C_3}Md\log n}{n}})^{\frac12} &\leq \frac{\kappa}{3}\max_{i \in [M]} \|\bmu^0_i-\bmu^*_i\|_2, \label{eq:ssize1}\\
({\textstyle\frac{\tilde{C_2}Md\log n}{n}})^{\frac12} &\leq \frac{\kappa}{12}.
\label{eq:ssize2}
\end{align}
The second inequality \eqref{eq:ssize2} in turn implies
\begin{align*}
&\textstyle\sup_{\bmu \in \cU} | \frac{1}{n} \sum_{j=1}^n w_i(X_j,\bmu) - \bbE\big[w_i(X;\bmu)\big]| \leq \frac{1}{12}\kappa
\end{align*}
for all $i\in[M]$. Let $\bmu^0$ be the initial iterate. We have 
\begin{align}
\|\bmu^1_i-\bmu^*_i\|_2 &= \frac{\|\frac{1}{n} \sum_{j=1}^nw_i(X_j;\bmu^0)(X_j-\bmu^*_i)\|_2}{\frac{1}{n} \sum_{j=1}^nw_i(X_j;\bmu^0)}\nonumber\\
&\leq \frac{\| \bbE w_i(X;\bmu^0)(X-\bmu^*_i) \|_2 + R_{\max} ({\frac{\tilde{C_3}Md\log n}{n}})^{\frac12}}{\bbE[w_i(X;\bmu^0)] - \frac{\kappa}{12}}\nonumber\\
&\overset{(i)}{\leq} \frac{\b\| \bbE w_i(X;\bmu^0)(X-\bmu^*_i) \b\|_2 + R_{\max}({\frac{\tilde{C_3}Md\log n}{n}})^{\frac12}}{ \frac{2\kappa}{3}} \nonumber\\
&\overset{(ii)}{\leq} \frac12 \max_{i \in [M]} \b\|\bmu^0_i-\bmu^*_i\b\|_2 + \frac{3R_{\max}(\tilde{C_3}Md\frac{\log n}{n})^{\frac12}}{2\kappa}, \label{eq:iterativecontraction}
\end{align}
where we appealed to
\begin{align*}
\bbE\big[w_i(X;\bmu^0)\big] &\geq \frac{3}{4}\kappa, \\
\| \bbE\big[w_i(X;\bmu^0)(X-\bmu^*_i)\big]\|_2 &\leq \frac{\kappa}{3} \max_{i \in [M]} \|\bmu^0_i-\bmu^*_i\|_2
\end{align*}
in steps (i) and (ii). Both are intermediate results established in the proof of Theorem \ref{cor:populationEMConvergence}. Finally, we have
\[
\begin{aligned}
\b\|\bmu^1_i-\bmu^*_i\b\|_2 &\leq \frac12 \max_{i \in [M]} \b\|\bmu^0_i-\bmu^*_i\b\|_2 +\frac{3R_{\max}(\tilde{C_3}Md\frac{\log n}{n})^{\frac12}}{2\kappa}\\
&\leq  \max_{i \in [M]} \b\|\bmu^0_i-\bmu^*_i\b\|_2,
\end{aligned}
\]
where the second step is a consequence of \eqref{eq:ssize1}. This implies $\bmu^1$ is also in the contraction region. By applying \eqref{eq:iterativecontraction} iteratively, we have
\begin{align}
\max_{i \in [M]} \|\bmu^t_i-\bmu^*_i\|_2 &\leq\frac{1}{2^t} \max_{i \in [M]} \|\bmu^0_i-\bmu^*_i\|_2\nonumber\\
&\qquad+\bb(1+\frac12 + \ldots + \frac{1}{2^{t-1}} \bb) \frac{3R_{\max}(\tilde{C_3}Md\frac{\log n}{n})^{\frac12}}{2\kappa} \\
&\leq 
\frac{1}{2^t} \max_{i \in [M]} \|\bmu^0_i-\bmu^*_i\|_2 + \frac{3R_{\max}(\tilde{C_3}Md\frac{\log n}{n})^{\frac12}}{\kappa}.
\end{align}

\section{Summary and discussion}

\textbf{Initialization.} We emphasize that our convergence results are local: they assume the EM algorithm is initialized in a neighborhood of the true centers. To obtain a such an initial iterate, we appeal to approaches based on the method of moments, such as the method proposed by \citet{Hsu13}. These methods are consistent; \ie\ as long as the sample size is large enough, they output a good initial iterate for the EM algorithm with high probability. 
Formally, as long as the sample size $n$ is greater than a polynomial of $d, M, 1/\epsilon, \log(1/\delta), 1/\kappa$, the method proposed by \citet{Hsu13} outputs estimates of the centers that satisfy
\[
\|\hat{\bmu}_i - \bmu^*_i \|_2 \leq (\|\bmu_i^*\|_2 + \lambda_{max})\epsilon \text{ with probability at least }1-\delta,
\]
where $\lambda_{max}$ is the largest eigenvalue of the covariance matrix of the samples.
The main benefit of combining a spectral method with the EM algorithm is the EM algorithm usually outputs a more accurate estimator. This is reflected in the smaller sample complexity of the EM algorithm.

\textbf{Minimum separation between centers.} Theorems \ref{thm:populationEMConvergence} and \ref{thm:EMConvergence} require the minimum separation between centers to grow as $\Omega(\min\{M,d\}^{\frac12})$. Compared to other methods for fitting mixtures of isotropic Gaussians, this dependence is suboptimal. For example, \citet{vempala04} showed that spectral clustering can accurately recover the labels in a mixture of spherical Gaussians provided that the minimum separation is at least $\Omega((M\log d)^{\frac14})$. Some approaches based on the method of moments are able to learn mixtures of Gaussians in which the centers are arbitrarily close together (as long as the sample size is large enough). However, the sample complexity of such methods are usually worse than that of the EM algorithm. 

If we restrict to studies of the EM algorithm and its variants (including gradient EM and the $K$-means algorithm), our requirement on the minimum separation between centers is optimal. \citet{bowei} imposes the same condition in their study of the convergence of the gradient EM algorithm. \citet{yulu} requires the minimum separation to grow proportionally to $M$ in their study of the convergence of the $K$-means algorithm. 

\textbf{Contraction radius and convergence rate.} This contraction radius in \eqref{eq:contractionradius} is optimal in the sense that it is approximately $\frac12 R_{min}$ when $R_{min}$ is large and we can find examples of the EM algorithm converging to non-global local minima if $\|\bmu_i-\bmu^*_i\|_2 = \frac12R_{\min}$ (see the remarks after Theorem \ref{thm:populationEMConvergence}). By comparison, the contraction radius for the gradient EM algorithm is very similar to \eqref{eq:contractionradius} \cite{bowei}, and the contraction radius for $K$-means is roughly $\frac12R_{\min} - CM^{\frac34}$ \cite[Theorem 6.2]{yulu}.

We show that the EM algorithm converges linearly up to statistical precision. This agrees with our simulation results and previous studies on the convergence of EM \citep{balakrishnan2017, tseng}. We see that the convergence rate decreases as $\kappa$ increases, which agrees with the folklore that the EM algorithm converges slowly on imbalanced mixture models. 

\textbf{Minimum sample size.} In terms of minimum sample size, our result is valid as long as $\frac{n}{\log n} \geq C\frac{MdR_{\max}^2}{\kappa^2 R_{\min}^2}$. \citet{bowei} established linear convergence of the gradient EM algorithm as long as $n \geq C\frac{M^6R_{\max}^6d}{R_{\min}^2}$ and \citet{yulu} established linear convergence of the $K$-mean algorithm as long as $\frac{n}{\log n} \geq C \frac{M}{\kappa^2}$. The variations in the minimum sample size are due to differences in the concentration results that appear in the proofs. We believe it is possible to avoid the $\log n$ factor and improve the dependence on $\frac{R_{\max}}{R_{\min}}$ in the minimum sample size by refining the concentration results in our proofs.

\appendix
\section{}
In appendix A, we prove lemma \ref{lemma:upperbound} and \ref{lemma:lowerbound}.

\subsection{Preliminaries}
Before jumping into the proof, we need the following preliminary result, which is lemma 12 and 13 in \citet{bowei}.
\begin{lemma}
  Suppose the minimum separation $R_{\min}$ and radius $a$ satisfy $R_{\min} \geq 30\sqrt{d}$ and 
  \[
  a \leq \frac12R_{\min} - \sqrt{d}\max\bb(4\sqrt{2[\log(R_{\min}/4)]_+},8\sqrt{3}\bb),
  \]
  then for any $\bmu$ such that $\bmu_i \in \cB(\bmu_i^*, a), \forall i \in [M]$, we have the following inequalities for any $p=0,1,2$ and any $i\neq j \in [M]$
  {\small 
    \begin{align*}
    \bbE w_i(X;\bmu)(1-w_i(X;\bmu))\|X-\bmu_i\|_2^p \leq 2M \bb( \frac32 R_{\max} + d \bb)^p 
    \text{\normalfont exp} \bbb( - \bb( \frac12R_{\min} -a \bb)  \frac{\sqrt{d}}{8}\bbb)\\
    \bbE w_i(X;\bmu)w_j(X;\bmu)\|X-\bmu_i\| \|X-\bmu_j\| \leq \frac{2-\kappa}{\kappa} \bb( \frac32 R_{\max} + d \bb)^2 
    \text{\normalfont exp} \bbb( - \bb( \frac12R_{\min} -a \bb)  \frac{\sqrt{d}}{8} \bbb).
    \end{align*}
  }
  \label{lemma:bowei}
\end{lemma}  
\subsection{Proof for lemma \ref{lemma:upperbound}}
We start from bounding $V_1$:
\begin{align}
V_1 \leq &\sup_{t \in [0,1]} \b\| \bbE w_1(X;\bmu^t)(1-w_1(X;\bmu^t))(X-\bmu^t_1)(X-\bmu^t_1)^T \b\|_{op} \nonumber \\
&+\sup_{t \in [0,1]} \b\| \bbE w_1(X;\bmu^t)(1-w_1(X;\bmu^t))(\bmu^*_1-\bmu^t_1)(X-\bmu^t_1)^T \b\|_{op} \nonumber\\
\leq&\sup_{t \in [0,1]} \b\| \bbE w_1(X;\bmu^t)(1-w_1(X;\bmu^t))(X-\bmu^t_1)(X-\bmu^t_1)^T \b\|_{op} \nonumber \\
&+a\sup_{t \in [0,1]} \b\| \bbE w_1(X;\bmu^t)(1-w_1(X;\bmu^t))(X-\bmu^t_1)^T \b\|_{2},\label{eq:V1bound}
\end{align}
where $a$ is the radius of the contraction region $\otimes_i \cB(\bmu^*_i,a)$. For any $t \in [0,1]$, there exists a rotation matrix $\Gamma$, such that all $\Gamma \bmu^t_i, \Gamma \bmu^*_i$, $i \in [M]$ have zero entries in the last $[d-2M]_+$ coordinates. Assume $d>2M$ for now, because this is the case where the rotation can yield a tighter bound. If $d \leq 2M$, this rotation is unhelpful but innocuous, and we can derive the same results without much modification.

Let $\tilde{X}=\Gamma X$, then $\tilde{X}|Z=i \sim \cN(\Gamma\bmu^*_i, I)$ and $\bbE \tX = \bbE X = 0$. Write
\[
\Gamma \bmu_i^t = [\tilde{\bmu}, \bzero_{d-2M}], \quad \Gamma \bmu^*_i = [\tilde{\bmu}^*_i, \bzero_{d-2M}], \quad \tilde{\bmu}^*_i,\tilde{\bmu} \in \bbR^{d-2M}. 
\]
It thus follows
{\small 
  \[
  (X-\bmu^t_1)(X-\bmu^t_1)^T=\Gamma^T \left[ 
  \begin{array}{cc}
  (\tX^{2M}-\tilde{\bmu}_1^t)(\tX^{2M}-\tilde{\bmu}_1^t)^T & (\tX^{2M}-\tilde{\bmu}_1^t)(\tX^{d-2M})^T\\
  (\tX^{d-2M})(\tX^{2M}-\tilde{\bmu}_1^t)^T & (\tX^{d-2M})(\tX^{d-2M})^T 
  \end{array}
  \right]\Gamma.
  \] 
}

Since $\tX^{d-2M}\sim \cN(0,I_{d-2M})$, it is independent of $\tX^M$. Also note that $w_1(X;\bmu^t)$ only depends on $\tX^{2M}$ (the part involving $X^{d-2M}$ cancels out), we have
\begin{align*}
&\b\| \bbE w_1(X;\bmu^t)(1-w_1(X;\bmu^t))(X-\bmu^t_1)(X-\bmu^t_1)^T \b\|_{op} = \left\| \left[ \begin{array}{cc}
D_1 & 0\\
0 & D_2
\end{array} \right] \right\|_{op} \\
\leq& \max(\|D_1\|_{op}, \|D_2\|_{op}).
\end{align*}
Applying lemma \ref{lemma:bowei} with dimension $\min\{2M,d\}$ and $p=2$, we have
\begin{equation}
\|D_1\|_{op} \leq 2M \bb( \frac32 R_{\max} + \min\{2M,d\} \bb)^2 \text{\normalfont exp} \bbb( - \bb( \frac12R_{\min} -a \bb)  \min\{d,2M\}^{\frac12}/8\bbb).
\label{eq:D1}
\end{equation}
Applying lemma \ref{lemma:bowei} with dimension $\min\{2M,d\}$ and $p=0$, we have
\begin{equation}
\begin{aligned}
\|D_2\|_{op} &= \bbE \bb[ w_1(\tX_{2M};\tilde{\bmu}^t)(1-w_1(\tX_{2M};\tilde{\bmu}^t)) \bb] \\
&\leq 2M \text{\normalfont exp} \bbb( - \bb( \frac12R_{\min} -a \bb)  \min\{d,2M\}^{\frac12}/8\bbb).
\label{eq:D2}
\end{aligned}
\end{equation}

Combining \eqref{eq:D1},\eqref{eq:D2} with \eqref{eq:V1bound}, we see
\begin{align}
V_1 &\leq 2M \bbb(\bb( \frac32 R_{\max} + \min\{2M,d\} \bb)^2 + a ( \frac32 R_{\max} + \min\{2M,d\} \bb) \bbb)  \nonumber\\
&\quad \cdot \text{\normalfont exp} \bbb( - \bb( \frac12R_{\min} -a \bb)  \min\{d,2M\}^{\frac12}/8\bbb) \nonumber \\
&\leq 2M \bb( 2 R_{\max} + \min\{2M,d\} \bb)^2 \text{\normalfont exp} \bbb( - \bb( \frac12R_{\min} -a \bb)  \min\{d,2M\}^{\frac12}/8\bbb).
\end{align}

Next we move to $V_i,i \neq 1$. Using the same decomposition
\begin{align}
V_i  \leq&\sup_{t \in [0,1]} \b\| \bbE w_1(X;\bmu^t)w_2(X;\bmu^t)(X-\bmu^t_1)(X-\bmu^t_i)^T \b\|_{op} \nonumber \\
&+a\sup_{t \in [0,1]} \b\| \bbE w_1(X;\bmu^t)w_i(X;\bmu^t)(X-\bmu^t_i) \b\|_{2}.
\end{align}
Apply the same rotation trick, we are able to show
\begin{equation}
V_i \leq \frac{2}{\kappa} \bb( 2 R_{\max} + \min\{2M,d\} \bb)^2 \text{\normalfont exp} \bbb( - \bb( \frac12R_{\min} -a \bb)  \min\{d,2M\}^{\frac12}/8\bbb).
\end{equation}
Since $\kappa \leq \frac{1}{M}$, we see
\[
\max_{i \in [M]} V_i \leq \frac{2}{\kappa} \bb( 2 R_{\max} + \min\{2M,d\} \bb)^2 \text{\normalfont exp} \bbb( - \bb( \frac12R_{\min} -a \bb)  \min\{d,2M\}^{\frac12}/8\bbb).
\]
The proof is now complete.

\subsection{Proof of lemma \ref{lemma:lowerbound}}
First, we can apply the same rotation trick to reduce the effective dimension to $\min\{M,d\}$. To do so, let $\Gamma$ be a rotation matrix such that the last $[d-M]_+$ coordinates of $\Gamma \bmu_i$ are zero for all $i \in [M]$. Write $\Gamma\bmu_i=(\tilde{\bmu}_i,\bzero_{[d-M]^+})$, $\tX = \Gamma X$ and we have
\[
\|X-\bmu_i\|_2^2=\|\Gamma X-\Gamma\bmu_i\|_2^2
=\|\tX_{M} - \tilde{\bmu}_i\|_2^2 + \|\tX_{[d-M]_+}\|_2^2.
\]
This implies $w_1(X;\bmu)=w_1(\tX_M;\tilde{\bmu})$ and $\bbE w_1(X;\bmu) = \bbE w_1(\tX_M;\tilde{\bmu})$ where $\tX_M|Z=i \sim \cN((\Gamma \bmu^*_i)_{M}, I_M)$. We thus have successfully reduced the effective dimension to $\min\{M,d\}$.

The rotation step is optional and can only reduce dimension when $d>M$. For ease of notation, let us assume $M \geq d$ and we opt not to do it. The next step in bounding $\bbE w_1(X;\bmu)$ is to restrict ourselves to the event where a) $X$ is generated by the first cluster, and b) $X$ lies in a ball $\cB(\bmu^*_1,r)$ for some radius to be selected later. Specifically, we have
\begin{equation}
\bbE w_1(X;\bmu)
\geq \pi_1 \bbE_{X \sim \cN(\bmu_1^*,I)} w_1(X;\bmu)
\geq \pi_1 \int_{\cB(\bmu_1^*, r)} w_1(x;\bmu)\phi(x;\bmu_1^*)dx.
\label{eq:sec6lbound1}
\end{equation}

Notice on $\cB(\bmu_i^*, r)$, by triangular inequality we have
\begin{eqnarray*}
  \|x-\mu_1 \|_2 &\leq& r+a\\
  \|x-\mu_i \|_2 &\geq& R_{\min}-r-a, \forall i \neq 1.
\end{eqnarray*}
Also, since $w_1(x;\mu)$ is decreasing in $\|x-\mu_1 \|_2$ and increasing in $\|x-\mu_i \|_2$, we have 
\begin{align*}
w_1(x;\mu)&\geq \frac{\pi_1 e^{-\frac{(r+a)^2}{2}}}{\pi_1 e^{-\frac{(r+a)^2}{2}} + (1-\pi_1) e^{-\frac{(R_{\min}-r-a)^2}{2}}}\\
&\geq 1 - \frac{1-\pi_1}{\pi_1} \exp{-\frac{1}{2}R_{\min}(R_{\min}-2r-2a)} \quad \\
&\geq 1 - \frac{1-\kappa}{\kappa} \exp{-\frac{1}{2}R_{\min}(R_{\min}-2r-2a)},
\end{align*}
where in the last step, we used numerical inequality $\frac{a}{a+b}\geq 1-\frac{b}{a}$. It thus follows
\begin{align}
& \quad \pi_1 \int_{\cB(\bmu_1^*, r)} w_1(x;\bmu)\phi(x;\bmu_1^*)dx \nonumber \\
& \geq \pi_1 \bigg(1-\frac{1-\kappa}{\kappa}\exp{-\frac{1}{2}R_{\min}(R_{\min}-2r-2a)}\bigg) \int_{\cB(\bmu_1^*, r)} \phi(x;\bmu_1^*)dx \nonumber \\
& \geq \kappa\bigg(1-\frac{1-\kappa}{\kappa}\exp{-\frac{1}{2}R_{\min}(R_{\min}-2r-2a)}\bigg) \bbP(\|\varepsilon \|_2 \leq r), \label{eq:sec6lbound2}
\end{align}
where $\varepsilon \sim \cN(\bzero, I_d)$. Moving forward. we naturally want to lower bound $\bbP(\|\varepsilon \|_2 \leq r)$, and the following lemma (lemma 8 in \cite{bowei}) allows us to do so. 

\begin{lemma}
  Let $X \sim \cN(\bzero, I_d)$, for $r \geq 2\sqrt{d}$, we have
  \[
  \bbP(\|X\|_2 \geq r) \leq \exp{-\frac{r\sqrt{d}}{2}}.
  \]
  \label{lemma:tail}
\end{lemma}

Let us pretend for now that $\frac{1}{2}R_{\min}(R_{\min}-2r-2a) \geq \frac{r\sqrt{d}}{2}$, then by chaining \eqref{eq:sec6lbound1}, \eqref{eq:sec6lbound2} and applying lemma \ref{lemma:tail}, we have
\begin{align*}
& \bbE w_1(X;\bmu)\\
\geq & \kappa\bigg(1-\frac{1-\kappa}{\kappa}\exp{-\frac{r\sqrt{d}}{2}}\bigg) \bigg(1-\exp{-\frac{r\sqrt{d}}{2}}\bigg) \\
\geq & \kappa - \exp{-\frac{r\sqrt{d}}{2}}.
\end{align*}
Therefore to let $\bbE w_1(X;\bmu) \geq \frac{3\kappa}{4}$, it suffices to let
\[
\exp{-\frac{r\sqrt{d}}{2}} \leq \frac{\kappa}{4}.
\]

Now we collect all the conditions we need for all the inequalities to go through; they are
\begin{itemize}
  \item[(C1)] $r \geq 2\sqrt{d}.$ (lemma  \ref{lemma:tail})
  \item[(C2)] $\frac{1}{2}R_{\min}(R_{\min}-2r-2a) \geq \frac{r\sqrt{d}}{2}$
  \item[(C3)] $\exp{-\frac{r}{2}\sqrt{d}} \leq \kappa/4$.
\end{itemize}
Setting $r = \frac{R_{\min}/2-a}{4}$, we can check
\begin{enumerate}
  \item $a \leq \frac{1}{2}R_{\min}-8\sqrt{d}$ implies (C1).
  \item $R_{\min} \geq \sqrt{d}/6$ implies (C2).
  \item $a \leq \frac{1}{2}R_{\min} - \frac{8}{\sqrt{d}}\log(\frac{4}{\kappa})$ implies (C3).
\end{enumerate}

If we replace all $d$ by $\min\{M,d\}$, the proof goes through with only notational changes. Now, we can readily check the conditions on $R_{\min}$ and $a$ in lemma \ref{lemma:lowerbound} imply the three conditions above. The proof is complete.

\section{}
In appendix B, we prove lemma \ref{lemma:contractionw} and lemma \ref{lemma:contractionvector} which facilitate the proof of Theorem \ref{thm:EMConvergence}. We first introduce some preliminary results on sub-gaussian random variables and then prove lemma \ref{lemma:contractionw} and lemma \ref{lemma:contractionvector}.

\subsection{Preliminaries}
\begin{lemma}
  Let $X$ be a random variable.
  \begin{enumerate}
    \item \textbf{(Sub-gaussian random variable).} $X$ is called sub-gaussian if there exists a finite $t$ such that $\bbE \exp{X^2/t^2} \leq 2$. For a sub-gaussian $X$, its sub-gaussian norm $\|X\|_{\psi_2}$ is defined as
    \[
    \|X\|_{\psi_2} = \inf\{t > 0: \; \bbE \exp{X^2/t^2} < \infty \}.
    \]
    
    \item \textbf{(Hoeffding's inequality).} Let $X_1,\ldots,X_N$ be independent, mean zero, sub-gaussian random variables. Then, for every $t\geq 0$, we have
    \[
    \bbP \bb(  \bb| \sum_{i=1}^N X_i \bb| \geq t \bb) \leq 2 \text{\normalfont exp}\bb( - \frac{ct^2}{\sum_{i=1}^N \|X_i\|_{\psi_2}^2} \bb),
    \]
    where $c$ is an absolute constant.
    \item \textbf{(Centering).} If $X$ is a sub-gaussian random variable, then $X-\bbE X$ is sub-gaussian too, and
    \[
    \|X - \bbE X\|_{\psi_2} \leq C\|X \|_{\psi_2},
    \]
    where $C$ is an absolute constant.
    \item \textbf{(Bounded random variable is sub-gaussian).} Any bounded random variable $X$ is sub-gaussian, with 
    \[
    \|X \|_{\psi_2} \leq C \|X\|_\infty,
    \]
    where $C=1/\sqrt{\log 2}$.
    \item Let $X$ be a bounded random variable on $[0,1]$, $Y$ be a sub-gaussian random variable. Then,
    \[
    \|XY\|_{\psi_2} \leq \|Y\|_{\psi_2}.
    \]
  \end{enumerate}
  \label{lemma:subgaussian}
\end{lemma}

\noindent \textbf{Proof of lemma \ref{lemma:subgaussian}:} Properties
1-4 are standard results from chapter 2 of \cite{hdp}; 5 is a consequence of taking expectation and infimum on both sides of
\[
\exp{\frac{X^2Y^2}{t^2}} \leq  \exp{\frac{Y^2}{t^2}}.
\]

\begin{lemma}
  Let $X$ be the mixture of $M$ unit variance gaussian distributions on $\bbR$, with centers denoted by $\{\theta_i\}_{i \in [M]}$ and mixing proportions by $\{\pi_i\}_{i \in [M]}$. Suppose $\max_{i \in [M]}|\theta_i|\leq R$ for some constant $R$. Then $X$ is sub-gaussian with sub-gaussian norm
  \[
  \|X\|_{\psi_2} \leq C \max(R,1) 
  \]
  for some absolute constant $C$.
  \label{lemma:subgaussianofGMM}
\end{lemma}

\noindent \textbf{Proof of lemma \ref{lemma:subgaussianofGMM}:} Let $Y$ be a random draw from centers $\{\theta_i\}_{i \in [M]}$ according to probabilities $\{\pi_i\}_{i \in [M]}$, i.e. $P(Y=\theta_i)=\pi_i$. It follows that $Y$ is a bounded random variable and $\|Y\|_{\psi_2} \leq C_1 R$. Let $\varepsilon \sim \cN(0,1)$. From standard results we know $\|\varepsilon\|_{\psi_2} \leq C_2$. Note that $X$ has the same distribution as $Y+\varepsilon$, we have 
\[
\|X\|_{\psi_2} = \|Y+\varepsilon\|_{\psi_2}
\leq \|Y\|_{\psi_2} + \|\varepsilon\|_{\psi_2}
\leq C_1 R + C_2 \leq C \max(R,1).
\]

\subsection{Proof of lemma \ref{lemma:contractionvector}}
Define $L:=1.5R_{\max}$, then for $a \leq \frac12R_{\min}$, we have $\otimes_i \cB(\bmu_i^*, a) \in \otimes_i \cB(0, L)$. This is a natural consequence of $\|\bmu^*_i\|_2 \leq R_{\max}$ for all $i$. To see why $\|\bmu^*_i\|_2 \leq R_{\max}$, suppose the opposite and, without loss of generality, let $\|\bmu^*_1\|_2 > R_{\max}$, then all $\bmu^*_i \in \cB(\bmu^*_1, R_{\max})$. Since $\bbE X = \sum \pi_i \bmu^*_i = 0$ but $\cB(\bmu^*_1, R_{\max})$ does not contain the origin, we get a contradiction. Also note that since theorem \ref{thm:populationEMConvergence} requires $R_{\min} \geq C_0\min\{d,2M\}^{\frac12}$ for a large $C_0$, we can work under the premise that $R_{\max} \geq 1$, because otherwise, even the population level convergence result does not apply.

Denote $\cU = \otimes_i \cB(0, L)$, and we establish all uniform convergence results on $\cU$. Let $n_\ve$ be the $\ve$-covering number of $\cB(0, L)$, then standard results \cite{hdp} have it $\log(n_\ve) \leq d \log(3L/\ve)$. By doing cartesian product on such covers, we can get a cover on $\cU$. We denote this cover by $M_\ve = \{\bmu^1,\ldots,\bmu^{N_\ve} \}$ with $M_\varepsilon \subset \bbR^{Md}$ and $\log(N_\ve) \leq Md \log(3L/\ve)$. For any $\bmu \in \otimes_i \cB(0, L)$, let $j(\bmu) = \argmin_{j \in [N_\ve]} \|\bmu - \bmu^j\|_2$. Then for all $\bmu \in \cU$, $\|\bmu_i - \bmu^{j(\bmu)}_i \|_2 \leq \ve$ for $\forall i \in [M]$.

Define
\[
I_{n,1}(\bmu) := \bb\| \frac{1}{n} \bb( \sum_{i=1}^nw_1(X_i;\bmu)(X_i-\bmu^*_1) \bb) - \bbE[w_1(X;\bmu)(X-\bmu^*_1)]\bb\|_2
\]
Start by noting
\begin{align*}
I_{n,1}(\bmu) \leq& \bb\| \frac{1}{n} \bb( \sum_{i=1}^nw_1(X_i;\bmu)(X_i-\bmu^*_1) \bb) - \frac{1}{n} \bb( \sum_{i=1}^nw_1(X_i;\bmu^{j(\bmu)})(X_i-\bmu^*_1) \bb) \bb\|_2 \\
&+ \bb\|\frac{1}{n} \bb( \sum_{i=1}^nw_1(X_i;\bmu^{j(\bmu)})(X_i-\bmu^*_1) \bb)-
\bbE[w_1(X;\bmu^{j(\bmu)})(X-\bmu^*_1)] \bb\|_2 \\
&+\bb\|
\bbE[w_1(X;\bmu^{j(\bmu)})(X-\bmu^*_1)] - \bbE[w_1(X;\bmu)(X-\bmu^*_1)] \bb\|_2 \\.
\end{align*}
Then we have
\[
\bbP \bbb( \sup_{\bmu \in \cU} I_{n,1}(\bmu) \geq t\bbb) \leq \bbP(A_t)
+\bbP(B_t)+\bbP(C_t),
\]
where the events $A_t,B_t,C_t$ are defined as
\begin{eqnarray*}
  &&A_t= \bbb\{ \sup_{\bmu \in \cU} \bb\| \frac{1}{n} \bb( \sum_{i=1}^nw_1(X_i;\bmu)(X_i-\bmu^*_1) \bb) - \frac{1}{n} \bb( \sum_{i=1}^nw_1(X_i;\bmu^{j(\bmu)})(X_i-\bmu^*_1) \bb) \bb\|_2 \geq \frac{t}{3} \bbb\},\\
  &&B_t= \bbb\{ \sup_{j \in [N_\ve]} \bb\|\frac{1}{n} \bb( \sum_{i=1}^nw_1(X_i;\bmu^{j})(X_i-\bmu^*_1) \bb)-
  \bbE[w_1(X;\bmu^{j})(X-\bmu^*_1)] \bb\|_2 \geq \frac{t}{3} \bbb\},\\
  &&C_t= \bbb\{ \sup_{\bmu \in \cU} \bb\|
  \bbE[w_1(X;\bmu^{j(\bmu)})(X-\bmu^*_1)] - \bbE[w_1(X;\bmu)(X-\bmu^*_1)] \bb\|_2 \geq \frac{t}{3} \bbb\}.\\
\end{eqnarray*}
For some $\delta > 0$, we next derive conditions on $t$ that suffice to let
\[
\bbP(A_t) \leq \frac{\delta}{2}, \quad \bbP(B_t) \leq \frac{\delta}{2}, \quad
\bbP(C_t) = 0.
\]
Then replacing $\delta$ with $\frac{1}{n}$ completes the proof.

\noindent\textbf{Upper bounding $\bbP(B_t)$}: 

Let $V_{1/2}$ be a $(1/2)$-cover of $\cB_d(0,1)$ with $\log|V_{1/2}|\leq d\log 6$. Then we know from standard result that
\begin{align*}
&\bb\|\frac{1}{n} \bb( \sum_{i=1}^nw_1(X_i;\bmu^{j})(X_i-\bmu^*_1) \bb)-
\bbE[w_1(X;\bmu^{j})(X-\bmu^*_1)] \bb\|_2\\
\leq& 2 \sup_{\bv \in V_{1/2}} \bb\langle\bv,\frac{1}{n} \bb( \sum_{i=1}^nw_1(X_i;\bmu^{j})(X_i-\bmu^*_1) \bb)-
\bbE[w_1(X;\bmu^{j})(X-\bmu^*_1)] \bb\rangle.
\end{align*}
Taking union bound, we have
\begin{align*}
&\quad \bbP(B_t)\\
 &\leq \bbP \bbb(  \sup_{\substack{j \in [N_\ve] \\ \bv \in V_{1/2}}} \bbb\{
\frac{1}{n} \sum_{i=1}^n 
\bb\langle\bv,  w_1(X_i;\bmu^{j})(X_i-\bmu^*_1) -
\bbE[w_1(X;\bmu^{j})(X-\bmu^*_1)] \bb\rangle
\bbb\}    \geq \frac{t}{6} \bbb)\\
&\leq \exp{Md \log(3L/\ve)+d\log 6}  \\
&\quad \cdot \sup_{\substack{j \in [N_\ve] \\ \bv \in V_{1/2}}} 
\bbP \bbb(   \bbb\{
\frac{1}{n} \sum_{i=1}^n 
w_1(X_i;\bmu^{j})\langle X_i-\bmu^*_1,\bv\rangle -
\bbE[w_1(X;\bmu^{j})\langle X-\bmu^*_1,\bv\rangle]
\bbb\}    \geq \frac{t}{6} \bbb).
\end{align*}
By part 5 of lemma \ref{lemma:subgaussian},
\[
\|w_1(X;\bmu^{j})\langle X-\bmu^*_1,\bv\rangle\|_{\psi_2} \leq \|\langle X-\bmu^*_1,\bv\rangle\|_{\psi_2}.
\]
Note that $\langle X-\bmu^*_1,\bv\rangle$ follows a one dimensional gaussian mixture model with centers at $\{\langle\bmu^*_i-\bmu^*_1,\bv\rangle \}_{i\in [M]}$. Since $|\langle\bmu^*_i-\bmu^*_1,\bv\rangle|\leq R_{\max}$ and $R_{\max} \geq 1$, we conclude from lemma \ref{lemma:subgaussianofGMM}
\[
\|\langle X-\bmu^*_1,\bv\rangle\|_{\psi_2} \leq CR_{\max} \leq CL.
\]
Consequently, by Hoeffding's inequality, we have
\[
\bbP(B_t) \leq \exp{Md \log(3L/\ve)+d\log 6-\frac{cnt^2}{L^2}}.
\]
To ensure $\bbP(B_t) \leq \frac{\delta}{2}$, it suffices to let
\[
t \geq C\sqrt{\frac{L^2\b(Md \log(\frac{18L}{\ve}) + \log(\frac2\delta)\b)}{n}}
\].

\noindent\textbf{Upper bound $\bbP(C_t)$}:

Using the same integration expression \eqref{eq:integrationexpression} as in section 4.1, we have
\begin{equation*}
\sup_{\bmu \in \cU} \bb\| 
\bbE[w_1(X;\bmu^{j(\bmu)})(X-\bmu^*_1)] - \bbE[w_1(X;\bmu)(X-\bmu^*_1)] \bb\|_2\\
\leq\ve \sum_{i=1}^{M} U_i,
\end{equation*}
where
\begin{eqnarray}
U_1 &=& \bbE \sup_{\bmu \in \cU} \b\|  w_1(X;\bmu)(1-w_1(X;\bmu))(X-\bmu_1^*)(X-\bmu_1)^T \b\|_{op},\\
U_i &=& \bbE \sup_{\bmu \in \cU} \b\|  w_1(X;\bmu)w_i(X;\bmu)(X-\bmu_1^*)(X-\bmu_i)^T \b\|_{op} \quad \text{   for } i \neq 1.
\end{eqnarray}
Since $C_t$ is deterministic, so as long as we have $\ve\sum_{i=1}^{M} U_i < \frac{t}{3}$, $C_t$ will never happen.

\noindent\textbf{Upper bound $\bbP(A_t)$:}

Using markov inequality, we have
\begin{align*}
&\quad \bbP(A_t) \\
&\leq \frac{3}{t} \bbE \bbb[ \sup_{\bmu \in \cU} \bb\| \frac{1}{n} \bb( \sum_{i=1}^nw_1(X_i;\bmu)(X_i-\bmu^*_1) \bb) - \frac{1}{n} \bb( \sum_{i=1}^nw_1(X_i;\bmu^{j(\bmu)})(X_i-\bmu^*_1) \bb) \bb\|_2 \bbb]\\
&\leq \frac{3}{t} \bbE \sup_{\bmu \in \cU} \bb\| 
\b( w_1(X;\bmu^{j(\bmu)}) - w_1(X;\bmu) \b) (X-\bmu^*_1)  \bb\|_2\\
&\leq \frac{3\ve}{t} \sum_{i=1}^{M} U_i. \quad \text{(due to mean value theorem)}
\end{align*}
To ensure $P(A_t) \leq \frac{\delta}{2}$, it suffices to let $t \geq 6\ve(\sum U_i)/\delta$. Note that whenever this holds, the condition that ensures $\bbP(C_t)=0$ is implied.

\noindent\textbf{Bounding $U_i$}:

Knowing that all $w \in [0,1]$, we see
\[
U_1 \leq \bbE \bb[ \sup_{\bmu \in \cU} \|X-\bmu_1^*\|\|X-\bmu_1\| \bb] \leq \bbE( \|X\| + L)^2 \leq \bbE( \|Y\| + \|\varepsilon \| + L)^2,
\]
where $Y$ is a random draw from centers $\{\bmu^*_i\}_{i \in [M]}$ according to probabilities $\{\pi_i\}_{i \in [M]}$ and $\varepsilon \sim \cN(\bzero, I_d)$.
With a bit more calculation, we see $U_1 \leq C'(L^2 + d)$, and the same bound also hold for other $U_i$. It thus follows
\[
\sum_{i=1}^{M} U_i \leq C'M(L^2 + d).
\]

\noindent\textbf{Conclusion:}

We set $\ve = \frac{\delta L}{6C'nM(L^2+d)}, \delta = \frac{1}{n}$, then any $t$ satisfy the following ensures bad events happen with probability less than $\delta$
\[
t \geq \max \bbb\{ \frac{L}{n}, CL\sqrt{\frac{Md \log \bb( \frac{108C'M(L^2+d)n}{\delta}  \bb) + \log \frac{2}{\delta}}{n}} \bbb\}.
\]
The second argument in maximum apparently dominates the first argument. After meticulously checking, we conclude  there exists a universal constant $C_3$, such that
\[
\bbP \bbb( I_{n,1}(\bmu) \geq R_{\max} \sqrt{\frac{\tilde{C_3}Md\log n}{n}}\bbb) \leq \frac{1}{n},
\]
where $\tilde{C_3}=C_3\log (M(3R_{\max}^2+d))$.

\subsection{Proof of lemma \ref{lemma:contractionw}}

The proof of lemma \ref{lemma:contractionw} is essentially the same as the proof of lemma \ref{lemma:contractionvector}. Start by noticing 
\begin{align*}
\bb|\frac{1}{n}\sum_{i=1}^nw_1(X_i;\bmu) - \bbE w_1(X;\bmu) \bb|
&\leq \bb| \frac{1}{n} \sum_{i=1}^nw_1(X_i;\bmu) - \frac{1}{n}  \sum_{i=1}^nw_1(X_i;\bmu^{j(\bmu)})\bb| \\
&+ \bb|\frac{1}{n} \sum_{i=1}^nw_1(X_i;\bmu^{j(\bmu)})-
\bbE w_1(X;\bmu^{j(\bmu)}) \bb| \\
&+\bb|
\bbE w_1(X;\bmu^{j(\bmu)}) - \bbE w_1(X;\bmu) \bb|.
\end{align*}
Then we have
\[
\bbP \bbb( \sup_{\bmu \in \cU} \bb|\frac{1}{n}\sum_{i=1}^nw_1(X_i;\bmu) - \bbE w_1(X;\bmu) \bb| \geq t\bbb) \leq \bbP(A_t)
+\bbP(B_t)+\bbP(C_t),
\]
where the events $A_t,B_t,C_t$ are defined as
\begin{eqnarray*}
  &&A_t= \bbb\{ \sup_{\bmu \in \cU} \bb| \frac{1}{n} \sum_{i=1}^nw_1(X_i;\bmu) - \frac{1}{n}  \sum_{i=1}^nw_1(X_i;\bmu^{j(\bmu)})\bb| \geq \frac{t}{3} \bbb\},\\
  &&B_t= \bbb\{ \sup_{j \in [N_\ve]} \bb|\frac{1}{n} \sum_{i=1}^nw_1(X_i;\bmu^{j(\bmu)})-
  \bbE w_1(X;\bmu^{j(\bmu)}) \bb| \geq \frac{t}{3} \bbb\},\\
  &&C_t= \bbb\{ \sup_{\bmu \in \cU} \bb|
  \bbE w_1(X;\bmu^{j(\bmu)}) - \bbE w_1(X;\bmu) \bb| \geq \frac{t}{3} \bbb\}.\\
\end{eqnarray*}
For some $\delta > 0$, we next derive conditions on $t$ that suffice to let
\[
\bbP(A_t) \leq \frac{\delta}{2}, \quad \bbP(B_t) \leq \frac{\delta}{2}, \quad
\bbP(C_t) = 0.
\]
Then replacing $\delta$ with $\frac{1}{n}$ completes the proof.

\noindent\textbf{Upper bounding $\bbP(B_t)$}: 

Since $w_1(X,\bmu)$ is bounded between [0, 1], it is sub-gaussian with $\|w_1(X,\bmu)\|_{\psi_2} \leq C$ for some absolute constant $C$. We can thus directly apply union bound and Hoeffding's inequality:
\begin{align*}
\bbP(B_t) &\leq \exp{Md \log(3L/\ve)} 
\cdot \sup_{j \in [N_\ve]} 
\bbP \bbb(   \bbb|
\frac{1}{n} \sum_{i=1}^n 
w_1(X_i;\bmu^{j}) -
\bbE w_1(X;\bmu^{j})
\bbb|    \geq \frac{t}{3} \bbb)\\
&\leq 2\exp{Md \log(3L/\ve)-cnt^2}.
\end{align*}
To ensure $\bbP(B_t) \leq \frac{\delta}{2}$, it suffices to let
\[
t \geq C\sqrt{\frac{Md \log(\frac{3L}{\ve}) + \log(\frac4\delta)}{n}}
\].

\noindent\textbf{Upper bound $\bbP(C_t)$}:

Using the same integration expression \eqref{eq:integrationexpression} as in section 4.1, we have
\begin{equation*}
\sup_{\bmu \in \cU} \bb| 
\bbE w_1(X;\bmu^{j(\bmu)}) - \bbE w_1(X;\bmu) \bb| 
\leq\ve \sum_{i=1}^{M} W_i,
\end{equation*}
where
\begin{eqnarray}
W_1 &=& \bbE \sup_{\bmu \in \cU} \b\|  w_1(X;\bmu)(1-w_1(X;\bmu))(X-\bmu_1) \b\|_{2},\\
W_i &=& \bbE \sup_{\bmu \in \cU} \b\|  w_1(X;\bmu)w_i(X;\bmu)(X-\bmu_i) \b\|_{2} \quad \text{   for } i \neq 1.
\end{eqnarray}
Since $C_t$ is deterministic, so as long as we have $\ve\sum_{i=1}^{M} U_i < \frac{t}{3}$, $C_t$ will never happen.

\noindent\textbf{Upper bound $\bbP(A_t)$:}

Using markov inequality, we have
\begin{align*}
\bbP(A_t) &\leq \frac{3}{t} \bbE \bbb[ \sup_{\bmu \in \cU} \bb| \frac{1}{n} \sum_{i=1}^nw_1(X_i;\bmu) - \frac{1}{n}  \sum_{i=1}^nw_1(X_i;\bmu^{j(\bmu)})\bb| \bbb]\\
&\leq \frac{3}{t} \bbE \bbb[ \sup_{\bmu \in \cU} \bb| 
w_1(X;\bmu^{j(\bmu)}) - w_1(X;\bmu)  \bb| \bbb]\\
&\leq \frac{3\ve}{t} \sum_{i=1}^{M} W_i. \quad \text{(due to mean value theorem)}
\end{align*}
To ensure $P(A_t) \leq \frac{\delta}{2}$, it suffices to let $t \geq 6\ve(\sum U_i)/\delta$. Note that whenever this holds, the condition that ensures $\bbP(C_t)=0$ is implied.

\noindent\textbf{Bounding $W_i$}:

Knowing that all $w \in [0,1]$, we see
\[
W_1 \leq \bbE \bb[ \sup_{\bmu \in \cU} \|X-\bmu_1\| \bb] \leq \bbE \b[ \|X\| + L \b] \leq \bbE \b[ \|Y\| + \|\varepsilon \| + L \b],
\]
where $Y$ is a random draw from centers $\{\bmu^*_i\}_{i \in [M]}$ according to probabilities $\{\pi_i\}_{i \in [M]}$ and $\varepsilon \sim \cN(\bzero, I_d)$.
With a bit more calculation, we see $U_1 \leq C'(L + \sqrt{d})$, and the same bound also hold for other $U_i$. It thus follows
\[
\sum_{i=1}^{M} U_i \leq C'M(L + \sqrt{d}).
\]

\noindent\textbf{Conclusion:}

We set $\ve = \frac{\delta }{6C'nM(L + \sqrt{d})}, \delta = \frac{1}{n}$, then any $t$ satisfy the following ensures bad events happen with probability less than $\delta$
\[
t \geq \max \bbb\{ \frac{1}{n}, C\sqrt{\frac{Md \log \bb( \frac{18C'ML(L + \sqrt{d})n}{\delta}  \bb) + \log \frac{4}{\delta}}{n}} \bbb\}.
\]
The second argument in maximum apparently dominates the first argument. After meticulously checking, we conclude  there exists a universal constant $C_2$, such that
\[
\bbP \bbb( \sup_{\bmu \in \cU} \bb| \sum_{i=1}^n w_1(X_i,\bmu) - \bbE w_1(X;\bmu) \bb| \geq  \sqrt{\frac{\tilde{C_2}Md\log n}{n}}\bbb) \leq \frac{1}{n},
\]
where $\tilde{C_2}=C_2\log (2MR_{\max} (2R_{\max}+\sqrt{d}))$.

\bibliographystyle{imsart-nameyear}
\bibliography{ruofei}

\end{document}